\documentclass[11pt,a4paper]{article}
\usepackage{amssymb}
\newcommand\qed{\hfill $\square$}

\newcommand{\NN}{\mathbb{N}}
\newcommand{\R}{\mathbb{R}}

\newcommand{\EE}{\mathcal{E}}

\newcommand{\DD}{\mathcal{D}}

\newtheorem{theo}{Theorem}

\newtheorem{lem}[theo]{Lemma}
\newtheorem{cor}[theo]{Corollary}
\newtheorem{rem}[theo]{Remark}

\newcommand{\beqn}{\begin{equation}}
\newcommand{\eeqn}{\end{equation}}
\newcommand{\bear}{\begin{eqnarray}}
\newcommand{\eear}{\end{eqnarray}}
\newcommand{\bean}{\begin{eqnarray*}}
\newcommand{\eean}{\end{eqnarray*}}

\newcommand{\Ac}{\mathcal{A}}
\newcommand{\Bc}{\mathcal{B}}
\newcommand{\Cc}{\mathcal{C}}
\newcommand{\Dc}{\mathcal{D}}

\newcommand{\OO}{\mathcal{O}}
\newcommand{\FF}{\mathcal{F}}

\newcommand{\eps}{\varepsilon}
\newcommand{\wto}{\rightharpoonup}

\begin{document}

\title{SELF-SIMILARITY FOR BALLISTIC AGGREGATION EQUATION}
\maketitle
\vspace{0.3cm}

\begin{center}
Miguel Escobedo\footnotemark[1],
St\'ephane Mischler\footnotemark[2]
\end{center}
\footnotetext[1]{Departamento de Matem\'aticas, Universidad del
Pa{\'\i}s Vasco, Apartado 644, E--48080 Bilbao, Spain.
E-mail~: {\tt mtpesmam@lg.ehu.es}}
\footnotetext[2]{Ceremade - UMR 7534,
Universit\'e Paris - Dauphine, Place du Mar\'echal De Lattre de 
Tassigny, 75775 Paris Cedex 16, France. 
    E-mail~: {\tt mischler@ceremade.dauphine.fr}}

\vspace{0.3cm}

\begin{abstract} We consider ballistic aggregation equation for gases in which each particle is identified either by its mass and impulsion 
or by its sole impulsion. For the constant aggregation rate we prove existence of self-similar solutions as well as convergence to the self-similarity for generic solutions. For some classes of mass and/or impulsion dependent rates we are also able to estimate the large time decay of some moments of generic solutions or to build some new classes of self-similar solutions. 
\end{abstract}

\vspace{0.3cm}

\section{Introduction}

\setcounter{equation}{0}
\setcounter{theo}{0}
In the present work, we are concerned with ballistic aggregation Smoluchowski like models for which we establish quantitative information on the qualitative behavior of solutions. By {\it ballistic aggregation}, also (improperly) called {\it kinetic coalescence} in previous works \cite{ELM,FM}, we mean that we consider a system of particles identified by their mass and impulsion which undergo an aggregation mechanism. That differs from the simplest aggregation mechanism introduced by Smoluchowski \cite{Sm2} in which model the particles are identified by their sole mass. 

\smallskip 
Let us be more precise. We denote by $P = P_{y}$ with $y = (m,p)$ a particle of mass $m > 0$ and impulsion $p \in \R^d$. The space of particles states is then $Y = \R_+ \times \R^d$ and the  velocity of the particle $P_y$ is $v=p/m$. We assume that at a microscopic level (the level of particles) the rate of collision of two particles $P = P_y$ and $P' = P_{y'}$ is a given nonnegative function $a = a(y,y')$ and when these two particles collide they join to form one aggregated particle $P'' = P_{y''}$ in such a way that the mechanism conserves total mass and total impulsion. In other words, the microscopic mechanism reads 
$$
P_{y} + P_{y'} \ \mathop{\longrightarrow}^{a(y,y')} \ P_{y''}, 
$$
with $y'' = (m'',p'')$ given by 
$$
m'' = m + m', \qquad p'' = p + p'.
$$
It is worth mentioning that the above reaction dissipates kinetic energy since that, denoting by $\EE^\sharp= m^\sharp \, |v^\sharp|^2/2$ the kinetic energy of particle $P^\sharp$,  we have
\bean
\EE^{**} - \EE - \EE^* &=& {1 \over 2} \, {|p+ p^*|^2 \over m+m_*}  -{1 \over 2} \, {|p|^2 \over m}  - {1 \over 2} \, {|p^*|^2 \over m_*} 
\\
&=& -  {1 \over 2} \, {m \, m_* \over m+m_*}  \, |v-v_*|^2. 
\eean

\smallskip
At the mesoscopic (or statistical or mean field) level, the system is described at time $t \ge 0$ by the density function $f(t,y) \ge 0$ of particles with state $y \in Y$. For a given initial distribution $f_{in}$, the evolution of  the density $f$ is described  by the Smoluchoswki/Boltzmann like equation:
\bear
\label{a1}
\partial_t f & = & Q(f)  \quad\hbox{ in }\quad
(0,+\infty)\times Y,\\
\label{a2}
f(0) & = &  f^{in} \quad\hbox{ in }\quad Y.
\eear
The collision operator $Q(f)$ is given by $Q(f) = Q_1(f) - Q_2(f)$, where
\bear
\label{q1}
Q_1(f)(y) & = & \frac{1}{2} \int_{\R^d}\int_0^m a(y',y-y')\ f(y')\
f(y-y')\ dm'dp',\\
\label{q2}
Q_2(f)(y) & = & \int_{\R^d}\int_0^\infty a(y,y')\ f(y)\ f(y')\ dm'dp'.
\eear

The two following examples of functions $a$ have been considered in relation with models in astrophysics \cite{W,BI}:
\beqn
\label{HS}
a(y,y') = a_{HS} (y,y') := (m^{1/3} + m'^{\, 1/3})^2 \, |v-v'|,
\eeqn
\beqn
\label{NP}
a(y,y') = a_{NP} (y,y') := {m+m' \over m m'} \, {1 \over |v-v'|^2}.
\eeqn
This model is seen as a simple test case or elementary analog of more realistic situations in fluid mechanics or astrophysics {\cite{CPY, JL}}. We refer to the introduction of \cite{RV,ELM,FM} for an elementary introduction to physics motivation of such a model. We also refer to \cite{CPY,JL,Trizac,W} and to the references quoted in \cite{RV,ELM,FM}  for a more detailed discussion about physics of aggregation.

In the context described above it is very natural to impose on the initial data $f_{in}$ to have finite number of particles and momentum. This condition reads:
\beqn
\label{a5}
0 \le f^{in}\in L^1\left( Y, (1+m+|p|) \, dydp
\right).
\eeqn
Existence of solutions under that condition has been proved in \cite{RV, ELM, FM}. It has also been proved that 
\beqn
\label{convL1}
f(t, \cdot) \to 0 \quad\hbox{in}\quad L^1 (Y),\,\,\,\hbox{as}\,\,\,t\to +\infty,
\eeqn
that is that the total number of particles tends to $0$. 

A more detailed description of the asymptotic behaviour of the solutions may be obtained by considering scaling invariance properties of the equations. This may be done for example by studying the so-called self similar solutions as it is possible to do for the Smoluchowski equation, see \cite{EM,FL}Ê
and  the references therein for recent results in that direction for that model. A first difficulty to this end is to determine the relevant scalings of the equation (\ref{a1}), (\ref{q1}), (\ref{q2}). We are very far from being able to treat the general case when the aggregation kernel $a(y, y')$ actually depends on both mass and momentum of the two colliding particles and even when the aggregation kernel $a(y, y')$ only depends on the momentum of the two colliding particles. We then may be less ambitious and just ask for whether a more accurate version than (\ref{convL1}) for some rate of aggregation $a$ is available? We may imagine to answer that question in several ways listed below by order of accuracy, and indeed depending of  the case we will establish any of such a kind of information.

\begin{itemize}
\item {\bf Answer 1. } 
Upper bound on moment: $\exists \, \bar\alpha$, $\exists \, \nu, \, C \in (0,\infty)$ such that 
$$
M_{\bar\alpha}(f(t,.)) \le {C \over t^\nu} \qquad \forall \, t \ge 1.
$$
\item {\bf Answer 2. } 
Upper and lower bound on moments: $\exists \,  \bar\alpha$, $\exists \, \nu_i = \nu_i( \bar\alpha), \, C_i= C_i( \bar\alpha)\in (0,\infty)$, such that 
$$
{C_1 \over t^{\nu_1}} \le M_{ \bar\alpha}(f(t,.)) \le {C_2 \over t^{\nu_2}} \qquad \forall \, t \ge 1.
$$
\item {\bf Answer 3. } 
Existence of self-similar solution:  there exists some profile function $\varphi_\infty : Y \to \R_+$, some exponents $\lambda, \, \mu, \, \nu \in \R$ such such that the function 
$$
\varphi(t,m,p) := t^\lambda \, \varphi_\infty(t^\mu \, m, t^\nu \, p)
$$
is a solution to equation (\ref{a1}), (\ref{q1}), (\ref{q2}).
\item {\bf Answer 4. } 
Convergence to self-similarity: for any given solution $f$ there exists a self-similar solution $\varphi$ such that $f \sim \varphi$ as $t \to \infty$, in a sense to be specified.
\end{itemize}

\noindent
Here depending of the model, we define the moment of order $\bar\alpha$ of $f$ in the following way:
\begin{itemize}
\item 
when $f = f(y)$ with $y = m \in Y = (0,\infty)$ or $y = p  \in Y = \R^d$, then $\bar\alpha = \alpha \in \R$ and
\beqn\label{def:Momenta}
M_{\bar\alpha}(f)  = M_{\alpha}(f) = \int_Y |y|^\alpha \, f \, dy;
\eeqn
\item 
when $f = f(y)$ with $y = (m,p) \in Y = (0,\infty) \times \R^d$, then $\bar\alpha = (\alpha,\beta) \in \R^2$ and
\beqn\label{def:Momentab}
M_{\bar\alpha}(f)  = M_{\alpha,\beta}(f) = \int_Y m^\alpha \, |p|^\beta \, f \, dy. 
\eeqn
\end{itemize}

\noindent
The results obtained in this work are very partials and may be classified as follows.

In Section 2 we consider the case of the kernel $a_{HS}(y, y')$ (which depends on both mass and momentum) and the only result we are able to prove is a upper estimate on some moments (that is a result of type "Answer 1"). 

In the remainder of the paper, we focus our attention on some toy models in which the aggregation rate $a$ depends upon the only impulsion or upon the only masses, namely $a(y,y') = a(p,p')$,  $a(y,y') = a(m,m')$ or even $a(y,y') \equiv 1$. The relation with the initial problem is not clear, and in particular it seems that a velocity depending aggregation rate $a(y,y') = a(v,v')$ should be more natural that an impulsion depending aggregation rate $a(y,y') = a(p,p')$. Anyway, on the one hand such kind of aggregation rates has been considered by physicists, see \cite{CPY,JL,Trizac,W}, and on the other hand our results and methods can give some ideas in order to tackle the so much more difficult models where the aggregation rate depends on both mass and momentum. 

Then, in Section 3 we consider a class of kernels which only depend on the momentum $p$ and $p'$, we establish some moment estimates of type "Answer 2", from which we deduce the rather strange conclusion that solutions do not enjoin a self-similar property (nor self-similar solution exists). That result sow doubt about the fact that in the case of the mass and impulsion hard spheres kernel,  solutions develop self-similar behavior. 

We treat in Section 4 the case where the kernel depends only on the masses $m$ $m'$ of the colliding particles and we exhibit a new class of self-similar solution (that is "Answer 3"). Lastly, in Section 5 the case of constant kernel is treated, for which results of type "Answer 3" and "Answer 4" are established. 

\smallskip
We end that introduction by some remarks and open questions. 
A common feature of these equations is that 
$$
M_{1,0}(t) \equiv M_{1,0}(0) \qquad \hbox{and}\qquad M_{0,0}(t) \to 0 \quad \hbox{as} \quad t \to \infty,
$$
and when the cross-section $a$ is homogeneous of order $\bar\gamma$ (which belongs to $\R$ or $\R^2$) it is likely that 
\beqn\label{Mgamma1t}
M_{\bar\gamma}(t) \, \equiv \, {1 \over t} \quad \hbox{as} \quad t \to \infty,
\eeqn
a result which is also known to be true for the coagulation equation (see \cite{FL,EMR,EM}) and for the inelastic Boltzmann equation (see \cite{MM3} and the references therein). The equivalence (\ref{Mgamma1t}) is established for the the impulsion depending and the mass depending aggregation rate, but only one side of that equivalence is proved in the case of the true mass and impulsion depending hard spheres aggregation rate. We ask then. 

\smallskip\noindent
{\bf Open question 1. }ÊIs it true that the assymptotic equivalence behavior (\ref{Mgamma1t}) holds for some true mass and impulsion depending aggregation rate? 

\smallskip\noindent
An other interesting question should be to establish some asymptotic behavior of  typical velocity or impulsion depending quantity. A way to express that in mathematical terms is the following: 

\smallskip\noindent
{\bf Open question 2. }ÊIs it possible to exhibit some moment $M_{\bar\alpha}$ for which we may determinate the long time behavior of  $M_{\bar\alpha}/M_0$ (even just saying that it converges to something)?


\section{Mass \& impulsion dependence case: a remark on the hard spheres model. }\label{sect:M&P}
\setcounter{equation}{0}
\setcounter{theo}{0}

Let us recall the following result
\begin{theo}\label{M&I:thHS} \cite[Theorem 2.6, Theorem 2.8 and Lemma 3.3]{FM} Assume that $a$ satisfies
\bean
&&0 \le a(y,y') = a(y',y) \le k_S(y) \, k_S(y'), \quad \forall \, y,y' \in Y,  
\\
&& a(m,-p,m',-p') = a(m,p,m',p') \quad \forall \, (m,p), (m',p')  \in Y,  
\\
&& a(m,p,m',p') \le a(m,p,m',-p') \quad \forall \, (m,p), (m',p')  \in Y \,\hbox{s.t.}\, \langle p,p' \rangle > 0,
\eean
with $k_S(y) := 1+m+|p|+|v|$.
For any even (in the $p$ variable) initial condition $0 \le f_{in}Ê\in L^1(Y; k_S^2(y) \,dy )$,  there exists a unique solution $f \in C([0,T); L^1(Y; k_S(y) \,dy )) \cap L^\infty (0,T; L^1(Y; k_S^2(y) \,dy ))$ $\forall \, T > 0$, which furthermore satisfies
\bear
\label{aHS:consmass}
&& \int_Y f (t,.) \, m \, dy \equiv \,\hbox{Cst}, \\
\label{aHS:consp}
&& f(t,.) \, \hbox{is even, so that } \,\, \int_Y f(t,.) \, p \, dy \equiv \,0, \\
\label{aHS:vkbdd}
&& \int_Y f (t,.) \, |v|^k \, dy \le \int_Y f_{in} \, |v|^k \, dy , \quad \forall \, k > 0, \\
\label{aHS:p2bdd}
&& \int_Y f (t,.) \,  |p|^2 \, dy \le \int_Y f_{in} \,  |p|^2 \, dy , \\
\label{aHS:Nbto0}
&& \int_Y f \, m^\alpha \, dy \,Ê\to \, 0 \quad\hbox{when}\quad t\to\infty, \quad \forall \, \alpha < 1. 
\eear
\end{theo}
 
 \begin{rem}\label{M&I:rkHS}
 (i) It is worth mentioning that the hard spheres collision rate $a_{HS}$ does satisfy the assumption of Theorem~\ref{M&I:thHS}, but not the Manev rate $a_{NP}$. 
 
 (ii) As a consequence of (\ref{aHS:consmass}), (\ref{aHS:vkbdd}), (\ref{aHS:p2bdd}) and (\ref{aHS:Nbto0}) we deduce that 
\beqn\label{def:Mab}
M_{\alpha,\beta}(t) := \int_Y f(t,.) \, m^\alpha \, |p|^\beta \, dy \,Ê\to \, 0 \quad\hbox{as}\quad t\to\infty
\eeqn
whenever $(\alpha,\beta)$ belongs to the region
$$
\{Ê\beta \in [0,2], \, \alpha < 1 -\beta/2 \} \cup \{ \beta \ge 2, \, \alpha < 2 -\beta  \}.
$$ 
 \end{rem}
 
In the case of the hard spheres model we are able to quantify the rate of decay of one of the moment functions of the solution. More precisely, we have the following result. 

\begin{lem}
\label{S2T1} Assume that $a = a_{HS}$. With the assumption of Theorem~\ref{M&I:thHS} there holds $A^{-1} := M_{-1/3,1}(0) < \infty$ and 
\beqn\label{eq:M-1/3,1}
M_{-1/3,1}(t) \le  \frac{1}{A+t/4} \qquad \forall \, t \ge 0.
\eeqn
\end{lem}

\noindent
\textbf{Proof of Lemma \ref{S2T1}.} First we have $M_{-1/3,1}(0) < \infty$ because
$$
m^{-1/3} \, |p| = m^{2/3} \, |v| \le m^{4/3} + |v|^2 \le 2 \, k_S^2.
$$
Now, from the expression (\ref{a1})-(\ref{a2}) of the collision kernel we have 
\bean
\int_Y Q(f,f) \, m^{-1/3} \, |p| \, dy =  {1 \over 2}Ê\int_Y \int_Y \Delta_{-1/3,\, 1} \, f \, f' \, dy dy',
\eean
with
\bean
\Delta_{-1/3,\, 1} = [Ê(m+m')^{-1/3} \, |p+p'| - m^{-1/3} \, |p| - (m')^{-1/3} \, |p'|] \, [r+r']^2 \, |v-v'|.
\eean
On one hand $-\Delta_{-1/3,\, 1} \ge 0$ because
$$
(m+m')^{1/3} \, \left( {|p| \over m^{1/3}} + {|p'| \over (m')^{1/3}} \right) \ge |p| + |p'|Ê\ge |p+p'|.
$$
On the other hand, if we only take into account the values of $v$ and $v'$ where $v\cdot v'<0$ and suppose that, for example, $|p|=\min(|p|, |p'|)$ we have

\bean
-\Delta_{-1/3,\, 1}&\ge& \left( { |p| \over m^{1/3}} + \left({ |p'| \over (m')^{1/3}} - { |p'| \over (m+m')^{1/3}}  \right)\right) \, [Êr^2 + (r')^2]Ê\, [|v| + |v'|] \\
&\ge& \left( { |p| \over m^{1/3}} \right) \, [(r')^2]Ê\, [ |v'|]  =  { |p| \over m^{1/3}} \, { |p'| \over (m')^{1/3}}.
\eean
Whence, by evenness of $f$ 
\bean
\label{S2estdelta}
{d \over dt} \int_Y f \,  { |p| \over m^{1/3}} \, dy 
&\le&- {1 \over 2} \int_{Y^2, \, v \cdot v' < 0}  { |p| \over m^{1/3}} \, { |p'| \over (m')^{1/3}} f \, f' \, dydy' \\
&\le&  - {1 \over 4}  \left(  \int_Y f \,  { |p| \over m^{1/3}} \, dy \right)^2,
\eean
from which (\ref{eq:M-1/3,1}) straightforwardly follows. \qed


\section{The  impulsion dependence case $a = a(p,p_*)$}

\setcounter{equation}{0}
\setcounter{theo}{0}

We consider now the equation (\ref{a1}), (\ref{q1}), (\ref{q2})  with a collision kernel $a$ independent of the mass of the colliding particles. We may then integrate the equation with respect to the mass and obtain that the function of $t$ and $p$, 
$
\int_0^\infty f(t, m, p)\, dm,
$
that we shall still denote $f$, satisfies the equation:
\bear
\label{S3a1}
\partial_t f & = & Q(f,f)  \quad\hbox{ in }\quad
(0,+\infty)\times \R^d,\\
\label{S3a2}
f(0) & = &  f_{in} \quad\hbox{ in }\quad \R^d, 
\eear
the collision operator $Q(f)$ is given by $Q(f,f) = Q_1(f,f) - Q_2(f,f)$, where
\bear
\label{S3q1}
Q_1(f,f)(y) & = & \frac{1}{2} \int_{\R^d} a(p',p-p')\ f(p')\ f(p-p') \, dp',\\
\label{S3q2}
Q_2(f,f)(p) & = & \int_{\R^d} a(p,p')\ f(p)\ f(p')\,dp'.
\eear
We focus on the cases 
\beqn\label{def:aggregp}
a(p,p') = |p- p'|^\gamma, \qquad \gamma \in [0,2], \quad d \in \NN^*.
\eeqn

Before stating our main result we need some definitions and notations.  We say that a function $f$ on $\R^d$ is even if 
$$
f(-p) = f(p) \qquad \forall \, p \in \R^d,
$$
it is radially symmetric if 
$$
f(R \, p) = f(p) \qquad \forall \, p \in \R^d, \,\, R \in SO(d)
$$
where $SO(d)$ stands for the rotation group on $\R^d$. 
For any weight function $k : \R^d \to \R_+$ we define the "moment of order $k$" of the non negative density measure$ f \in M^1_{loc}(\R^d)$ by 
$$
M_k (f) := \int_{\R^d}  k(p) \, f(dp),
$$
and we define $M^1_k$ as the set of Radon measures $\mu$ such that $M_k(|\mu|) < \infty$. For any $\alpha \in \R_+$ we use the shorthand
notation 
$$
M_\alpha :=\int_{\R} f(p)\, |p|^\alpha\, dp,
$$
that is $M_\alpha = M_k(f)$ for $k(p) = |p|^\alpha$ and the shorthand notation $M^1_\alpha = M^1_\ell$ for $\ell(p) = 1 + |p|^\alpha$.

\begin{theo}\label{theo:aggregp} Consider the aggregation rate (\ref{def:aggregp}).

\smallskip
(i) For any even initial datum $f_{in}Ê\in M^1_{2\alpha}$, $\alpha \in \NN \backslash \{ 0, 1Ê\}$, there exists a unique even solution $f \in C([0,T);M^1(\R^d)-weak) \cap L^\infty(0,T;M^1_{2\alpha}(\R^d))$ to equation (\ref{S3a1})--(\ref{S3q2}).  For any $\alpha \in [0,1]$ the function $t \mapsto M_\alpha(t)$ is decreasing and $f(t,.)$ is radially symmetric for any $t \ge 0$ if furthermore $f_{in}$ is radially symmetric.

\smallskip
(ii) Moreover, the solution $f(t,.)$ satisfies 
\beqn\label{estim:Mgammap}
{1 \over M_\gamma(0)^{-1} + k_1 \, t} \le M_\gamma (t) \le {1 \over M_\gamma(0)^{-1} + k_2 \, t} \qquad \forall \, t \ge 0,
\eeqn
for some constants $k_i = k_i(\gamma,d) \in (0,\infty)$. 
\end{theo}
One of the main tools in order to establish that result is to consider moment equations. As it is classical for the coagulation equation, but here using one more change of variable $p' \to -p'$, any even solution $f$ to equation (\ref{S3a1})--(\ref{S3q2}) satisfies (at least formally) the fundamental moment equation
\bear
{d \over dt} M_\alpha &=&
{1 \over 2}Ê \int_{\R^d} \!  \int_{\R^d} f \, f' \, a(p,p') \, [|p+p'|^\alpha - |p|^\alpha - |p'|^\alpha]Ê\, dpdp' \nonumber \\
&=&
{1 \over 4}Ê \int_{\R^d} \!  \int_{\R^d} f \, f' \, \left\{ a(p,p') \, \left[|p+p'|^\alpha - |p|^\alpha - |p'|^\alpha \right]Ê\right. \nonumber\\
&&\hskip 2cm \left.
+ \,a(p,-p') \, \left[|p-p'|^\alpha - |p|^\alpha - |p'|^\alpha \right]Ê\right\} \, dpdp'.  \label{S3EAka}
\eear

More precisely, we consider in this Section the  case $\gamma \in (0,2)$ and $d \in \NN^*$, the case $\gamma = 1$ and $d = 1$ and the case $\gamma=2$ and $d \in \NN^*$. The case $\gamma=0$ and $d=1$ is treated in Section~5. We shall use the following notation for the moments of order $\alpha \in \NN$:

\subsection{Proof of the existence and uniqueness part in Theorem~\ref{theo:aggregp}.}

We prove in this subsection an uniqueness and existence result for a general class of aggregation rates by adapting some arguments from \cite{LM,FM}, see also \cite{Norris}. We then deduce the existence and uniqueness part in Theorem~\ref{theo:aggregp}.

\begin{lem}
\label{lem:gama2Uniq} We consider a continuous aggregation rate $a : \R^{2d}Ê\to \R_+$ which satisfies 
\bear\label{eq:struct-a1}
a(-p,-p') &=& a(p,p') \qquad \forall \, p,p' \in \R^d,
\\ \label{eq:struct-a2}
a(p,p') &\le& a(-p,p') \qquad \forall \, p,p' \in \R^d, \,\, p \cdot p' > 0,
\eear
a even weight function $k : \R^d \to \R_+$ and we define 
$$
\Delta_k (p,p') := a (p,p') \, [k(p'') + k(p') - k(p)], \quad \tilde \Delta_k (p,p') =A (p,p') + A (-p,p') .
$$
We assume that 
\beqn\label{lem:gama2Uniq2}
a(p,p') \le C \, k(p) \, k(p') 
\quad \hbox{and} \quad \tilde A (p,p') \le  C \, k(p) \, k(p')^2.
\eeqn
For any given even initial datum $f_{in} \in M^1_k(\R^d)$ there exists no more than one even solution $f \in C([0,T);M^1_k(\R^d)) \cap L^\infty(0,T;M^1_{k^2}(\R^d))$ to equations  (\ref{S3a1})--(\ref{S3q2}). \end{lem}

\begin{rem}\label{rem:gama2Uniq} (i) The same result holds without the evenness assumption on the density function when the second condition in (\ref{lem:gama2Uniq2}) is replaced by 
$$
A (p,p') \le  C \, k(p) \, k(p')^2.
$$
We refer to  \cite{LM,FM} where such kind of result is proved in a $L^1$ framework. 
The same result also holds for radially symmetric solutions when we assume that 
\beqn\label{eq:struct-aRadial}
a(Rp,Rp') = a(p,p') \quad \forall \, p,p' \in \R^d, \,\, R \in SO(d),
\eeqn
and  the second condition in (\ref{lem:gama2Uniq2}) is replaced by 
$$
\int_{R \in SO(d)} A (p,R \, p') dR \le  C \, k(p) \, k(p')^2.
$$
(ii) The same kind of result holds for aggregation rate defined on $Y^2$ with $Y = (0,\infty) \times \R^d$ as it is the case when particles are identified by their mass and impulsion, see \cite{FM}.
\end{rem}

\noindent\textbf{Proof of Lemma \ref{lem:gama2Uniq}.} {\sl Step 1. }ÊWe claim that for $g \in C([0,T); M^1_{k}-weak)$, $G \in L^1(0,T; M^1_{k})$ and  $b \in C((0,T)Ê\times \R^d; \R_+)$  such that 
\beqn\label{eq:dtg=G-ag}
\partial_t g = G - b \, g \quad\hbox{in the sense of } \,\, \DD'([0,T) \times \R^d),
\eeqn
the differential inequality 
\beqn\label{ineq:dtg=G-ag}
{d \over dt} \|Êg \, k \|_{M^1}  \le \|  G \, k \|_{M^1}Ê- \|Êb \, g \, k \|_{M^1}Ê
\eeqn
holds in the sense of $\DD'([0,T))$. First, it is clear using a classical duality argument that equation (\ref{eq:dtg=G-ag}) has at most one solution. Indeed, given two solutions $g_1, g_2 \in  C([0,T); M^1_{k}-weak)$, we have for any $t \in (0,T)$, $\varphi_t \in C_{comp}(\R^d)$ and denoting by $\varphi \in C_{comp}([0,t] \times \R^d)$ the solution to the dual homogeneous equation $\partial_t \varphi = b \, \varphi$ 
$$
\int_{\R^d} (g_2 - g_1)(t) \, \varphi_t \, dp = \int_0^t\int_{\R^d} \{Ê(\partial_s g_2 - \partial_s g_1) \, \varphi +( g_2 -  g_1) \, \partial_s \varphi \}Ê\, dsdp = 0.
$$
Now, for any $g_\eps (0) \in C_{K_\eps} := \{Êu \in C(\R^d); \,\, \hbox{supp} \, u \subset K_\eps \}$, with $K_\eps \subset \R^d$ a compact, and any 
$G_\eps \in L^1(0,T;C_{K_\eps})$ there exists a (unique) solution $g_\eps \in C([0,T);C_{K_\eps})$ to equation (\ref{eq:dtg=G-ag}) which furthermore satisfies 
\bear\nonumber
{d \over dt} \int_{\R^d} |g_\eps|Ê\, k \, dy 
&=& \int_{\R^d} (G_\eps - b \, g_\eps) \, \hbox{sign} g_\eps \, k \, dy \\ \label{ineq:dtg=G-ag2}
&\le& \int_{\R^d} |G_\eps|  \, k \, dy - \int_{\R^d} |g_\eps|Ê\, a \, k \, dy.
\eear
Here, $\hbox{sign} g_\eps = 1$ if $g_\eps > 0$, $\hbox{sign} g_\eps =  0$ if $g_\eps = 0$, $\hbox{sign} g_\eps = -1$ if $g_\eps < 0$. Finally, we can build (by a standard truncation and regularization by convolution process) the sequences $(G_\eps)$ and $g_\eps(0)$ such that furthermore $G_\eps \wto G$, $g_\eps(0) 
\wto g(0)$ in the weak sense of measures in $M^1_k$, $\|ÊG_\eps (s) \|_{M^1_k} \le \|ÊG(s) \|_{M^1_k}$ for a.e. $s \in (0,T)$, $\|Êg_\eps (0) \|_{M^1_k} \le \|Êg(0) \|_{M^1_k}$. By the previous uniqueness argument we have $g_\eps \wto g$ in the weak sense of measure and we get (\ref{ineq:dtg=G-ag}) by passing to the limit in (\ref{ineq:dtg=G-ag2}). 

\smallskip\noindent
 {\sl Step 2. } Let us consider two solutions $f_1, f_2 \in C([0,T);M^1_k(\R^d)) \cap L^\infty(0,T;M^1_{k^2}(\R^d))$ which are evens and let us denote $D = f_2 - f_1$, $S= f_1 + f_2$.  By a standard algebraic computation $D$ satisfies the following equation
\bean
\partial_t D 
&=& \hat Q(f_2,f_2)  - \hat Q(f_1,f_1)  = \hat Q(D,S) \\
&=& \hat Q_1(D,S) - S \, L(D) -  L(S) \, D,
\eean
where 
$$
 \hat Q_i (\varphi,\psi) = {1 \over 2}Ê\, ( Q_i (\varphi,\psi) + Q_i (\psi,\varphi)), 
 \quad L(\varphi) := \int_{\R^d} a(p,p') \, \varphi(p') \, dp'.
$$
Because of the assumption made on $a$ and $f$ we have $D \in C([0,T];M^1_k-weak)$, $G := \hat Q_1(D,S) - S \, L(D) \in L^\infty(0,T;M^1_{k})$ and $0 \le b := L(S) \in C([0,T]Ê\times \R^d)$ so that the first step implies
\bean
{d \over dt} \| D \|_{M^1_k}Ê
&\le& \| (\hat Q_1(D,S) - S \, L(D)) \, k \|_{M^1}Ê- \|ÊD \, k \, L(S)Ê\|_{M^1}
\\
&\le& {1 \over 2} \int \!\! \int a \, [k'' + k'] \, |D(dp)|Ê\, S(dp') -  {1 \over 2} \int \!\! \int a \, k \, |D(dp)|Ê\, S(dp') 
\\
&\le&{1 \over 4}  \int \!\! \int \tilde A \, |D(dp)|Ê\, S(dp') \le  {C \over 4} \, \|ÊS \|_{M^1_{k^2}} \,  \|ÊD \|_{M^1_{k}}.
\eean
Uniqueness follows by using the Gronwall lemma. \qed

\begin{lem}
\label{lem:gama2Exist} Consider a continuous aggregation rate $a : \R^{2d}Ê\to \R_+$ which satisfies (\ref{eq:struct-a1}) (resp. (\ref{eq:struct-aRadial})), (\ref{eq:struct-a2}) as well as 
\beqn\label{eq;structa3}
a(p,p') \le C \, (k + k') \qquad \forall \, p, \, p' \in \R^d,
\eeqn
for the weight function $k(p) = 1 + |p|^2$ and some constant $C \in (0,\infty)$. 
For any given even  (resp. radially symmetric) initial datum $f_{in} \in M^1_{2\alpha}(\R^d)$ there exists at least one even (resp. radially symmetric) solution $f \in C([0,T);M^1(\R^d)-weak) \cap L^\infty(0,T;M^1_{2\alpha}(\R^d))$ to equation (\ref{S3a1})--(\ref{S3q2}), and this one furthermore satisfies  $t \mapsto M_\beta (t)$ is decreasing for any $\beta \in [0,1]$.
 \end{lem}

\begin{rem}
\label{rem:gama2Exist} It is likely that by adapting some arguments introduced in \cite{MW}, see also \cite{EM0,LM}, for any even  (resp. radially symmetric) initial datum $f_{in} \in L^1_{2\alpha}(\R^d)$ the approximating solution $f_n(t,.)$ built in the proof below is a Cauchy sequence in $C([0,T;L^1(\R^d))$ so that we may conclude $f \in C([0,T);L^1(\R^d)) \cap L^\infty(0,T;L^1_{2\alpha}(\R^d))$. 
\end{rem}

\noindent\textbf{Proof of Lemma \ref{lem:gama2Exist}.} We define the sequence of bounded aggregation rates $a_n := a \wedge n$, for which classically fixed point argument (see for instance \cite{FM} which deals with some similar situation) implies the existence of a unique even  (resp. radially symmetric) solution $f_n \in C([0,T);L^1_{2\alpha}(\R^d))$ to equation (\ref{S3a1})--(\ref{S3q2}) associated with $a_n$ for any initial datum $f_{in,n} \in L^1_{2\alpha+2}(\R^d)$, $\alpha \in \NN$, $\alpha \ge 2$. 
Then, we have for any $\beta \in \NN^*$, $\beta \le \alpha$
\bean
{d \over dt} \int f_n \,  (1 + |p|^{2\beta})
&=& {1 \over 2}Ê\int f_n \, f_n' \, a_n \,   \left[ (|p|^2 + 2 \, p \cdot p' + |p'|^2)^\beta - |p|^{2\beta} - |p'|^{2\beta} - 1 \right]Ê \\
&=& Ê\int f_n \, f_n' \, a_n \,   \left[ 2 \, \beta \, p \cdot p' \,  |p|^{2 (\beta-1)} -1/2 \right]Ê\\
&&+ \sum \mu_{\beta_1,\beta_2,\beta_2} \int f_n \, f_n' \, a_n \,   (p \cdot p')^{\beta_1} \,  |p|^{2 \, \beta_2} \, \,  |p|^{2 \, \beta_3},
\eean
where in the last sum we have $\beta_1 + \beta_2 + \beta_3 = \beta$ and ($\beta_1 \ge 2$ or ($\beta_2 \ge 1$ and $\beta_3 \ge 1$)) or, in other words, 
$ |p \cdot p'|^{\beta_1} \,  |p|^{2 \, \beta_2} \, \,  |p|^{2 \, \beta_3} \le |p|^{2 \, \beta'} \,  |p'|^{2 \, ( \beta-\beta')}$ with $1 \le \beta' \le \beta-1$. Since we also have
\bean
&& Ê\int f_n \, f_n' \, a_n \, p \cdot p' \,  |p|^{2 (\beta-1)}  = \\
&&\qquad =  Ê\int_{p \cdot p' > 0}  f_n \, f_n' \, (a(p,p') \wedge n - a(-p,p') \wedge n) \, p \cdot p' \,  |p|^{2 (\beta-1)} \le 0,
\eean
we conclude with 
\beqn\label{ineq:fnLbeta}
{d \over dt} \int f_n \,  (1 + |p|^{2\beta})
\le \sum_{1 \le \beta' \le \beta - 1}Ê \mu_{\beta'} \int f_n \, f_n' \, a \,    |p|^{2 \, \beta'} \,  |p'|^{2 \, ( \beta-\beta')}.
\eeqn
When $\beta = 1$ the set of admissible values of $\beta'$ is empty, and we recover a result from \cite{EM}
$$
{d \over dt} \int f_n \,  (1 + |p|^{2}) \le 0,
$$
so that 
\beqn\label{ineq:fnLk}
\sup_{[0,T]} \|Êf_n \|_{L^1_k} \le  \|Êf_{in,n} \|_{L^1_k}.
\eeqn
When $\beta \ge 2$, gathering (\ref{eq;structa3}), (\ref{ineq:fnLbeta}) and  (\ref{ineq:fnLk}), we easily conclude by a iterative argument that 
\beqn\label{ineq:fnLkbeta} 
\sup_{[0,T]} \|Êf_n \|_{L^1_{k^\beta}} \le C_T(\beta,  \|Êf_{in,n} \|_{L^1_{k^\beta}}).
\eeqn
Considering a sequence $(f_{in,n})$ such that  $f_{in,n} \wto f_{in}$ in the weak sense of measure and $\|Êf_{in,n} \|_{L^1_{k^\beta}}$ remains bounded, we easily pass to the limit in the equation satisfied by $f_n$ thanks to (\ref{ineq:fnLkbeta}). 
The fact that $t \mapsto M_\beta(t)$ is decreasing comes from the fact that $p \mapsto |p|^\beta$ is a sub-additive function when $\beta \in [0,1]$, so that $\Delta_\beta \le 0$ and then $d/dt M_\beta(t) \le 0$. \qed

\medskip
\noindent\textbf{Proof of the existence and uniqueness part in Theorem~\ref{theo:aggregp}. } It is clear that $a(p,p') = |p-p'|^\gamma$ satisfies (\ref{eq:struct-a1}), (\ref{eq:struct-a2}), the first inequality in (\ref{lem:gama2Uniq2}) and (\ref{eq;structa3}). Moreover, the second inequality in (\ref{lem:gama2Uniq2}) holds since we have 
\bean
\tilde \Delta_2 
&= & |p-p_*|^\gamma \, ( |p+p_*|^2 + |p_*|^2 - |p|^2 + 1) +  |p+p_*|^\gamma \, ( |p-p_*|^2 + |p_*|^2 - |p|^2 + 1) 
\\
&= & 2 \,  ( |p-p_*|^\gamma -  |p+p_*|^\gamma ) \, p \cdot p'  + ( |p-p_*|^\gamma +  |p+p_*|^\gamma ) \, ( 2 \, |p_*|^2  + 1),
\eean
where the first term in non positive and the second term is bounded by say $8 \, (k')^2 \, k$, using that $|p\pm p_*|^\gamma \le 2 \, (|p|^\gamma + |p'|^\gamma)$. 
We conclude by using Lemma~\ref{lem:gama2Uniq} and Lemma~\ref{lem:gama2Exist}.  \qed

\subsection{Proof of the rate decay part in Theorem~\ref{theo:aggregp} when $\gamma < 2$.}

For an even initial datum $f_{in} \in M^1_4(\R^d)$ we consider the unique even solution $f \in C([0,T);M^1-weak) \cap L^\infty(0,T;M^1_4)$, $\forall \, T$, given by Theorem~\ref{theo:aggregp}Ê(i). This one satisfies the moment equation
\beqn\label{eq:Mgamma}
{d \over dt} M_\gamma =
{1 \over 2}Ê \int_{\R^d} \!  \int_{\R^d} f \, f' \, \Delta_\gamma \, dpdp' = {1 \over 4}Ê \int_{\R^d} \!  \int_{\R^d} f \, f' \, \tilde\Delta_\gamma \, dpdp',
\eeqn
with
$$
- \Delta_\gamma = |p-p'|^\gamma \, [|p+p'|^\gamma - |p|^\gamma - |p'|^\gamma]Ê
$$
and 
\beqn\label{def:tDeltaGamma}
- \tilde\Delta_\gamma =  |p-p'|^\gamma \, \left[|p+p'|^\gamma - |p|^\gamma - |p'|^\gamma \right]Ê
+ |p+p'|^\gamma \,  \left[|p-p'|^\gamma - |p|^\gamma - |p'|^\gamma \right]Ê . 
\eeqn
We split the proof of  Theorem~\ref{theo:aggregp}Ê(ii) in several steps. 

\smallskip\noindent
{\sl Step 1. } One the one hand, for any given $A > 0$ and any $p,p' \in \R^d$ such that $A^{-1}Ê\, |p'| \le |p|Ê\le A \, |p|$
we easily get
\bear\nonumber
|\Delta_\gamma|Ê
&\le& (|p|+|p'|)^\gamma \, \maxÊ\left[ (|p|+|p'|)^\gamma, |p|^\gamma + |p'|^\gamma \right] \\ \label{estim:region1}
&\le& 2^4 \, \max(|p|,|p'|)^{2\gamma} \le 2^4 \, A^\gamma \, (|p| \, |p'|)^\gamma.
\eear
On the other hand, we define $M := \max(|p|,|p'|)$, $m := \min(|p|,|p'|)$, $x := m/M \in [0,1]$, $\eps := \hat p \cdot \hat p' \in [-1,1]$ and we compute (in the first line we have assumed that $|p|= M$ which is not a restriction to the generality because of the symmetry of $\tilde\Delta_\gamma$)
\bear\nonumber
-\tilde\Delta_\gamma 
&=& M^{2\gamma} \, \left\{   |\hat p-x \, \hat p'|^\gamma \, \left[1  + x^\gamma - |\hat p+ x \, \hat p'|^\gamma \right]Ê
+  |\hat p+x \, \hat p'|^\gamma \, \left[1  + x^\gamma - |\hat p- x \, \hat p'|^\gamma \right]Ê Ê\right\} \\ \nonumber
&=& M^{2\gamma} \, \left\{ (1+x^\gamma) \, [ (1 + 2 \, \eps \, x + x^2 )^{\gamma/2} + (1 - 2 \, \eps \, x + x^2 )^{\gamma/2} ]Ê \right. \\ \nonumber
&&\qquad\qquad\left. - 2 \, (1 + 2 \, \eps \, x + x^2 )^{\gamma/2} \,  (1 - 2 \, \eps \, x + x^2 )^{\gamma/2} \right\} \\ \label{estim:region2}
&=& M^{2\gamma} \, \left\{ 2 \, x^\gamma + \OO (x^2) \right\}  \le 3 \,  M^{2\gamma} \, x^\gamma = 3 \, (|p| \, |p'|)^\gamma 
\eear
uniformly on $\eps \in [-1,1]$ and $x \le A_0^{-1}$ for $A_0 \ge 1$ large enough. 

\smallskip
Gathering (\ref{estim:region1}) and (\ref{estim:region2}) we obtain
$$
{1 \over 4} \, Ê\tilde \Delta_\gamma \ge -  k_1 \, |p|^\gamma \, |p'|^\gamma \qquad \forall \, p,p' \in \R^d,
$$
with $k_1 :=  \max(3/4,2^3 \, A_0^\gamma)/4 $, and equation (\ref{eq:Mgamma}) then  implies
$$
{d \over dt} M_\gamma  \ge - k_1 \, M_\gamma^2.
$$
We straightforwardly obtain the first inequality in (\ref{estim:Mgammap}) by integrating this differential equation.

\smallskip\noindent
{\sl Step 2. } First, together with the variables $M$, $x$ and $\eps$ introduced in Step 1, we define $r > 0$ and $u \in [0,1]$ by setting $r^2 := |p|^2 + |p'|^2$ and $u := 2 \, p \cdot p'/r^2$, so that 
$|p\pm p'|^2 =  r^2 \, (1 \pm u)$. Splitting the positive and the negative terms in identity (\ref{def:tDeltaGamma}), we have 
\bean
-\tilde\Delta_\gamma 
&=& (|p|^\gamma + |p'|^\gamma)  \, (  |p-p'|^\gamma + |p+p'|^\gamma ) -  2 \, |p-p'|^\gamma \,  |p+p'|^\gamma 
\\
&=& r^{2\gamma} \, \left\{Ê{ (|p|^2)^{\gamma/2} + (|p'|^2)^{\gamma/2} \over  (|p|^2 + |p'|^2)^{\gamma/2}} \, \left[Ê( 1 + u )^{\gamma/2} +  (1 - u )^{\gamma/2} \right] 
- 2 \, ( 1 + u )^{\gamma/2} \,   (1 - u )^{\gamma/2}  \right\} .
\eean
Since $\gamma/2 \in [0,1]$, the map $x \mapsto x^{\gamma/2}$ is sub-additive, and we obtain 
\bean
-\tilde\Delta_\gamma 
&\ge & r^{2\gamma} \, \left\{Ê \left[Ê( 1 + u )^{\gamma/2} +  (1 - u )^{\gamma/2} \right]  - 2 \, ( 1 + u )^{\gamma/2} \,   (1 - u )^{\gamma/2}  \right\} 
\\
&\ge & M^{2\gamma} \,  ( 1 + u )^{\gamma/2} \,   (1 - u )^{\gamma/2}  \, \phi(u), \quad
\phi(u) :=  \left[Ê( 1 - u )^{-\gamma/2} +  (1 + u )^{-\gamma/2} \right]  - 2.
\eean
We easily verify that $\phi$ is increasing on $[0,1]$ so that $\phi(u) > \phi(0) = 0$ for any $u \in [-1,1]$, $u \not = 0$. Coming back to the variables $M$, $x$ and $\eps$,  that is $\phi(u) >  0$ for any $p, p' \in \R^d$ such that the associated variables $M$, $x$ and $\eps$ satisfy $M > 0$, $ x > 0$ and $\eps \not= 0$. 
Moreover, when $\eps = 0$ ($p$ and $p'$ are orthogonal vectors) we also have 
\bean
-\tilde\Delta_\gamma
&=& 2 (|p|^2 + |p'|^2)^{\gamma/2} \, \left[ |p|^\gamma + |p'|^\gamma -  (|p|^2 + |p'|^2)^{\gamma/2} \right]  \\
&\ge& 2 \, M^{2\gamma} \, \left[ 1 + x^\gamma -  (1 + x^2)^{\gamma/2} \right] > 0
\eean
for any $p,p' \in \R^d$ such that the associated variables $M$ and $x$ satisfy $M > 0$, $ x > 0$, because the function $z \mapsto z^{\gamma/2}$ is strictly sub-additive on $\R_+$, that is $(z+z')^{\gamma/2} < z^{\gamma/2} + (z')^{\gamma/2}$ for any $z,z' > 0$.  Gathering these two lower bounds on $-\tilde\Delta_\gamma$, it yields
\beqn\label{estim:LowerBddDeltaGamma1}
-\tilde\Delta_\gamma \ge M^{2 \, \gamma} \, \psi(x,\eps)
\eeqn
with $ \psi(x,\eps) > 0$ for any $x > 0$ and $\eps \in [-1,1]$. 

\smallskip
Next, coming back to (\ref{estim:region2}), we also deduce 
\beqn\label{estim:LowerBddDeltaGamma2}
-\tilde\Delta_\gamma 
= M^{2\gamma} \, \left\{ 2 \, x^\gamma + \OO (x^2) \right\}  \ge  M^{2\gamma} \, x^\gamma
\eeqn
uniformly on $\eps \in [-1,1]$ and $x \le A_0^{-1}$ for $A_0 \ge 1$ large enough. Gathering (\ref{estim:LowerBddDeltaGamma1}) with (\ref{estim:LowerBddDeltaGamma2}) we deduce that for some constant $k_2 > 0$ we have 
$$
\forallÊ\, p, p' \in \R^d \qquad -{1 \over 4}Ê\tilde\Delta_\gamma  \ge  k_2 \,  M^{2\gamma} \, x^\gamma = k_2 \,  (|p| \, |p'|)^\gamma,
$$
and  equation (\ref{eq:Mgamma}) then  implies
$$
{d \over dt} M_\gamma  \le - k_2 \, M_\gamma^2.
$$
The second inequality in (\ref{estim:Mgammap}) is again obtained by integrating this differential equation.

\subsection{The case $a (y,y') = |p-p'|$, $d=1$.}

In the particular case under consideration $d=1$ and $\gamma=1$, we can establish a more accurate version of the decay estimate on 
the moment $M_1$ together with additional moment estimates.

\begin{lem}
\label{S3TA0} Assume $a(y, y') = |p-p'|$ and $d=1$. For any even initial datum $f_{in}Ê\in M^1_3(\R)$ the unique solution $f \in C([0,T];M^1(\R)) \cap L^\infty(0,T; M^1_3(\R))$ of (\ref{S3a1})-(\ref{S3q2}) given by Theorem~\ref{theo:aggregp} satisfies for any $t \ge 0$
\bear
\label{S3M0}
\max \left( \frac{M_0(0)}{(1+M_1(0)\, t/2)^2},\frac{2^{3/2}\,M_0(0)}{(2+3M_3^{1/3}(0)\, t)^{3/2}}  \right) 
&\le& M_0(t) \,\, \le \,\, \frac{M_0(0)}{(1+M_1(0)t)^{1/2}} \qquad \qquad \\
\label{S3M1}
{1 \over M_1(0)^{-1} + t} &\le&  M_1(t) \,\, \le  \,\, {1 \over M_1(0)^{-1} + t/2} \\
\label{S3M2estbis}
{ M_2(0)  \over (1 + M_1(0) \,  t/2)^{2}} &\le& M_2(t) \,\,\le\,\, M_2(0) \\
\label{S3M3est}
 { M_3(0)  \over (1 + M_1(0) \,  t/2)^{2}} &\le& M_3(t) \,\,\le\,\, M_3(0). 
\eear
\end{lem}

\begin{rem}
\label{S3Rem2}
The above estimates on the behaviour of $M_1(t)$ for $t$ large are quite good. That is not the case for the estimates on $M_\alpha$, $\alpha = 0, 2, 3$ which seem to be rather partial. Worst, with these bounds we can not even know what is the limit of any of the quotients of moments $M_\alpha(t)/M_1(t)$ for $\alpha = 0, 2, 3$ as $t\to \infty$. The value of such a limit would indicate whether the solution $f(t)$ has a tendency to concentrate or to spread as $t$ increases (see also below the discussion concerning the case $\gamma = 2$).
\end{rem}

\noindent{\bf Proof of Lemma~\ref{S3TA0}. }ÊIntroducing the notations $M = \max(|p|,|p'|)$, $m = \min (|p|, |p'|)$, we systematically exploit the differential equation
\bear
{d \over dt} M_\alpha 
=  {1 \over 4}Ê \int_{\R} \!  \int_{\R} f \, f' \, \Delta_\alpha \, dpdp'
\eear
with 
$$
\Delta_\alpha := [M-m]Ê\, [ (M+m)^\alpha - M^\alpha - m^\alpha] + [M+m]Ê\, [ (M-m)^\alpha - M^\alpha - m^\alpha].
$$

\smallskip\noindent{\sl Step 1. $\alpha = 1$. }ÊWe have 
$$
\Delta_1 = - 2 \, (M+m) \, m,
$$
from which we deduce
$$
\frac{d}{dt}M_1(t)=-\frac{ M_1^2(t)}{2}-\frac{B_1(t)}{2}, \qquad
B_1(t) := \int_\R \!  \int_{\R} f \, f' \, \{\min(|p|,|p'|)Ê\} ^2\, dpdp'. 
$$ 
Since $0 \le \{\min(|p|,|p'|)Ê\} ^2 \le |p|Ê\, |p'|$, we have $0 \le B_1(t) \le M^2_1(t)$ and we obtain the two closed differential inequalities
\bear
 -M_1^2(t) \le \frac{d}{dt}M_1(t)\le -\frac{M_1^2(t)}{2},
\eear
from which we deduce (\ref{S3M1}). 

\smallskip\noindent{\sl Step 2. $\alpha = 0$. }ÊWe have 
$$
\Delta_0 = - 2 \, M,
$$
from which we deduce
\beqn\label{ed:ddtM0}
\frac{d}{dt}M_0(t)=- {B_0(t) \over 2}, \qquad B_0(t) := Ê\int_\R \!  \int_{\R} f \, f' \,\max(|p|,|p'|) \, dpdp'. 
\eeqn 
Since $|p|Ê\le \max(|p|,|p'|) \le |p| + |p'|$, we have $M_0\, M_1\le B_0 \le 2\,M_0\, M_1$ and then
\bear
\label{S3M0M1est}
-ÊM_0 \, M_1 \le {d \over dt}ÊM_0 \le -Ê{1 \over 2}Ê\, M_0 \, M_1.
\eear
Using the previous estimate (\ref{S3M1}) on $M_1(t)$ we get 
\bean
 -\frac{M_0(t)}{M_1^{-1}(0)+t/2}  \le \frac{d}{dt}M_0(t)\le -\frac{M_0(t)}{2(M_1^{-1}(0)+t)},
\eean
from which we deduce the first lower estimate as well as the upper bound in (\ref{S3M0}).

\smallskip\noindent{\sl Step 3. $\alpha = 2$. }ÊÊWe have 
$$
\Delta_2 = - 4 \, m \, M^2
$$
from which we deduce
$$
\frac{d}{dt}M_2(t)=- B_2(t), \qquad
B_2(t) := \int_\R \!  \int_{\R} f \, f' \, \min (|p|,|p'|)Ê \, |p|Ê\, |p'| \, dpdp'. 
$$
Using that $0\le \min (|p|,|p'|)Ê \, |p|Ê\, |p'| \le |p|^2 \, |p'|$ together with  (\ref{S3M1}), we obtain
$$
- M_2 \,  {1 \over M_1(0)^{-1} + t/2} \le - M_2 \, M_1 \le \frac{d}{dt}M_2(t) \le 0,
$$
which implies (\ref{S3M2estbis}).

\smallskip\noindent{\sl Step 4. $\alpha = 3$. }ÊÊÊWe have 
$$
0 \ge \Delta_3 = - 2 \, M \, m^3 - 2 \, m^4 \ge - 4 \, M \, m^3 \ge - 4 \, |p|^3 \, |p'|,
$$
from what we deduce
$$
0 \ge {d \over dt}ÊM_3(t) \ge - M_1 \, M_3,
$$
which again implies (\ref{S3M3est}). 

\smallskip\noindent{\sl Step 5. $\alpha = 0$ again.}Ê
Coming back to the moment $M_0$, we write for any $\eps > 0$
 \bean
\frac{d}{dt}  M_0
&=&
- {1 \over 2}Ê\int_\R \!  \int_{\R} f \, f' \, | p' - p | \, dpdp' \\
&\ge&
- {1 \over 4}Ê\int_\R \!  \int_{\R} f \, f' \, (\eps + {1 \over \eps} \, |p-p'|^2) \, dpdp' \\
&\ge&
- {\eps \over 4} \, M_0^2 - {2 \over \eps} \, M_0 \, M_2 . 
\eean
By interpolation we have $M_2(t) \le M_0^{1/3}(t)Ê\, M_3^{2/3}(t)$. Since, by (\ref{S3M3est}), $M_3(t)\le M_3(0)$ for all $t>0$Ê 
we deduce $M_2(t) \le M_0^{1/3}(t)Ê\, M_3^{2/3}(0)$. Therefore 
\bean
\frac{d}{dt}M_0(t)
&\ge&
- {\eps \over 4} \, M_0^2 - {2 \over \eps} \, M^{4/3}_0 \, M_3^{2/3}(0)
\eean
We now chose $\eps\equiv\eps(t) > 0$ such that 
${\eps} \, M_0^2 = {1\over \eps} \, M^{4/3}_0 \, M_3^{2/3}(0)$, or equivalently $ {\eps}  =  M^{-1/3}_0\, M_3^{1/3}(0)$.
With that choice of $\eps(t)$ the moment equation reads 
\bean
\frac{d}{dt}M_0(t)
\ge  - \frac{9}{4}\,  M_3^{1/3}(0)\, M^{5/3}_0,
\eean
from which we deduce the second lower estimate in (\ref{S3M0}).
\qed

\begin{rem}
\label{step5} In the last step, we may also argue as follows. Gathering the estimate $\max(|p|,|p'|) \ge (|p| \, |p'|)^{1/2}$, the differential equation (\ref{ed:ddtM0})
and the interpolation estimate $M_1^{5/2} \le M_{1/2}^2 \, M_3^{1/2}$ we obtain thanks to (\ref{S3M3est}) 
$$
{d \over dt}ÊM_0 \le - {1 \over dt}Ê\,  M_3^{1/2}(0) \, M_1^{5/2}(t).
$$
Together with (\ref{S3M1}) we recover the second lower estimate in (\ref{S3M0}).
\end{rem}


\subsection{The case $a = |p-p_*|^2$}

In the particular case under consideration $\gamma=2$ and $d \in \NN^*$, we can close the family of moment equations for any moments $M_{2\alpha}$, $\alpha \in \NN$. In the following lemma we give the expression of moments up to order $4$, showing a (unexpected?) non self-similar behavior of solutions.

\begin{lem}
\label{S3gamma2} Assume $a(y, y') = |p-p'|^2$and $d \in \NN^*$. There exists a numerical constant $k_d \in (0,\infty)$, $k_1 := 2$,  such that for any radially symmetric initial datum $f_{in}Ê\in M^1_6(\R)$ the unique radially symmetric solution $f \in C([0,T];M^1(\R)) \cap L^\infty(0,T; M^1_6(\R))$ of (\ref{S3a1})-(\ref{S3q2}) given by Theorem~\ref{theo:aggregp} satisfies for any $t \ge 0$
\bear
\label{S4M0}
M_0 (t) &=&  {M_0(0) \over (M_2(0)^{-1} + 2 \, k_d \,  t)^{1/(2k_d)}} \\
\label{S4M2}
M_2 (t) &=& {1 \over M_2(0)^{-1} + 2 \, k_d \,  t} \\
\label{S4M4}
M_4 (t) &=& M_4(0) \,  ( M_2(0)^{-1} + 2 \, k_d \,  t )^{1/k_d - 2} . 
\eear
\end{lem}

\noindent\textbf{Proof of Lemma \ref{S3gamma2}.} We proceed in several steps. 

\smallskip\noindent{\sl Step 1. $\alpha = 2$. }ÊUsing the fact that $f$ is radially symmetric (so that the odd moments of $f$ vanish) and the notations $p = r \, \sigma$, $r = |p|$, $p' = r' \, \sigma'$, $r' = |p'|$, the fundamental moment identity (\ref{S3EAka}) implies 
\bean
{d \over dt}ÊM_2 
&=& {1 \over 2} \int_{\R^d} \!\! \int_{\R^d} f \, f' \, [|p|^2 - 2 \, p \cdot p' + |p'|^2] \, (2 \, p \cdot p') \, dpdp' \\
&=& - 2  \int_{\R^d} \!\! \int_{\R^d} f \, f' \, [p \cdot p' ]^2  \, dpdp' \\
&=& - 2  \int_0^\infty \int_0^\infty f(r) \, f(r') \, r^{d+1}Ê\, (r')^{d+1} \, dr dr' Ê\times \int_{S^{d-1}} \!\! \int_{S^{d-1}}[\sigma \cdot \sigma' ]^2  \, d\sigma d\sigma' \\
&=& - 2 \, k_d \, M_2^2,
\eean
with 
\bean
k_d &:=& \left(  \int_{S^{d-1}} \!\! \int_{S^{d-1}}[\sigma \cdot \sigma' ]^2  \, d\sigma d\sigma' \right) \times \hbox{meas}(S^{d-1})^{-2} \\
&=&  \hbox{meas}(S^{d-1})^{-1} \int_{S^{d-1}}\sigma_1^2   \, d\sigma . 
\eean
We compute $k_1 = 1$, $k_2 = 1/2$. 
The expression (\ref{S4M2}) immediately follows by integrating that ODE. 

\smallskip\noindent{\sl Step 2. $\alpha = 0$. }ÊWhen $\alpha =0$, the fundamental moment identity (\ref{S3EAka}) and the fact that $f$ is radially symmetric imply 
\bean
{d \over dt}ÊM_0 
&=& {1 \over 2} \int_{\R^d} \!\! \int_{\R^d} f \, f' \, [|p|^2 - 2 \, p \cdot p' + |p'|^2] \, (-1) \, dpdp' \\
&=& - M_2 \, M_0.
\eean
Integrating that ODE with the help of (\ref{S4M2}) we get (\ref{S4M0}). 

\smallskip\noindent{\sl Step 3. $\alpha = 4$. }ÊWhen $\alpha = 4$, the fundamental moment identity (\ref{S3EAka}) and the fact that $f$ is radially symmetric imply
\bean
{d \over dt}ÊM_4 
&=& {1 \over 2} \int_{\R^d} \!\! \int_{\R^d} f \, f' \, [|p|^2 - 2 \, p \cdot p' + |p'|^2] \, [Ê4 \, (p \cdot p')^2 + 8 \, |p|^2 \, (p \cdot p') + 2 \, |p|^2 \, |p'|^2]Ê\, dpdp' \\
&=& {1 \over 2} \int_{\R^d} \!\! \int_{\R^d} f \, f' \, \left\{Ê[2 \, |p|^2] \, [Ê4 \, (p \cdot p')^2 + 2 \, |p|^2 \, |p'|^2]Ê- 16 \, |p|^2 \, (p \cdot p')^2
\right\} \, dpdp' \\
&=& 2 \int_{\R^d} \!\! \int_{\R^d} f \, f' \, \left\{Ê |p|^4 \, |p'|^2 - 2 \, |p|^2 \, (p \cdot p')^2
\right\} \, dpdp' \\
&=& (2 - 4 \, k_d) \, M_2 \, M_4.
\eean
Integrating that ODE with the help of (\ref{S4M2}) we get (\ref{S4M4}). 
\qed

\begin{rem}
\label{S3Rem3} (i) On the one hand, the moment $M_\alpha(g(t,.))$ of a self-similar function $g$ of the form $g(t,p) = t^\mu \, G(t^\nu \, p)$ satisfies 
$$
M_\alpha(g(t,.)) = C_\alpha \, t^{\mu - (d+\alpha) \, \nu}.
$$
On the other hand,  when $d=1$ we have $k_1 = 1$ so that the solution $f$ of equation (\ref{S3a1})-(\ref{S3q2})  satisfies
$$
M_0(f(t,.)) \sim C_0' \, t^{-1/2}, \quad M_2(f(t,.)) \sim C_2' \, t^{-1}, \quad  M_4(f(t,.)) \sim C_4' \, t^{-1}.
$$
Since the long time behavior of these functions are incompatibles, there does not exist any self-similar solution with self-similar profile $G \in M^1_6(\R)$.  

\noindent
(ii) When $d=1$, to make the ideas simpler, the moment $M_{2\alpha}$ satisfies the edo 
$$
{d \over dt} M_{2\alpha} = \sum_{\beta = 1}^{\alpha-1}Ê\pmatrix{2\alpha \cr 2\beta} \, M_{2\beta} \, M_{2(\alpha+1-\beta)} - 
\sum_{\beta = 0}^{\alpha-1}Ê\pmatrix{2\alpha \cr 2\beta+1} \, M_{2\beta+2} \, M_{2(\alpha-\beta)}.
$$
In particular, we find 
$$
{d \over dt}ÊM_6 = 3 \, M_2 \, M_6 - 5 \, M_4^2.
$$
When $M_2(0) = 1/2$ (for the sake of simplification again), the solution is
$$
M_6(t) = \left( M_6(0) - 2 \, M_4(0)^2 + {2 \, M_4(0)^2 \over (1+t)^{5/4}} \right) \, (1+t)^{3/2} \qquad \forall \, t \ge 0,
$$
with $M_6(0) - 2 \, M_4(0)^2  > 0$ (Holder inequality). The solutions of equation (\ref{S3a1})-(\ref{S3q2})  have a rather strange behavior
since that 
$$
M_0 \sim \kappa_0 \, t^{-1/2}, \quad {M_2 \over M_0}  \sim \kappa_1 \, t^{-1/2}, \quad  {M_6Ê\over M_0}  \sim \kappa_2 \, t^{3/2},
$$
In some sense, the behavior is in part comparable with the solutions of the inelastic Boltzmann equation which energy (here the $M_2$ moment) dissipates and  in part comparable with the solutions to Smoluchoski equation which high moments rapidly increase. It is worth mentioning that here the "mean second moment" (that is $M_2/M_0$) tends to $0$ in the large time assymptotic. The opposite feature occurs for some models dealt in section~\ref{sect:mass}. 
\end{rem}


\section{The mass dependence case $a = a(m,m_*)$}
\label{sect:mass}

\setcounter{equation}{0}
\setcounter{theo}{0}
Consider now the problem (\ref{a1})-(\ref{q2}) where the kernel $a(y, y')$ only depends on the masses of the particles, namely
\bear
\label{S4E1}
a(y, y')= a(m, m'),
\eear
and introduce the associated Smoluchowski equation 
\bear
&&\frac{\partial F}{\partial t}(t, m)=
\frac{1}{2}\int_0^m\, F(t, m-m') F(t, m')\, a(m-m', m')\, dm' \nonumber\\
&&\hskip 5cm -\int_0^\infty  F(t, m) F(t, m')\, a(m, m')\, dm'. 
\label{S4E2}
\eear
For any function $\psi \in L^1(\R^3)$ we define the Fourier transform $\FF$ and the inverse Fourier transform $\FF^{-1}$ by 
$$
\hat \psi (\eta) = (\FF \psi) (\eta) = \int_{\R^3} \psi(p) \, e^{- i \, p \cdot \eta} \, dp, 
\qquad
(\FF^{-1} \psi) (p) = (2 \, \pi)^{-3} \int_{\R^3} \psi(\eta) \, e^{ i \, p \cdot \eta} \, dp.
$$

\begin{theo}
\label{S4T1}
For any smooth function $a$ on $\R^3$ homogeneous of degree $\theta^{-1}$, $\theta \in (0,\infty)$, and such that $\varphi := \FF^{-1} (e^{-a(\cdot)}) \ge 0$, and for any solution $F\equiv F(t, m)$ to the coagulation equation (\ref{S4E2}) with coagulation kernel $a(m, m')$, the function $f(t, m, p)$ defined by
\bear
&&f(t, m, p)=m^{-3\, \theta}\, F(t, m)\,  \varphi\left(\frac{p}{m^\theta}\right) \label{S4E1bis},
 \label{S4E1ter}
\eear
is a solution of the equation (\ref{a1}), (\ref{q1}), (\ref{q2}) for the same aggregation kernel.
\end{theo}

\begin{rem}
\label{S4R1} Theorem \ref{S4T1} is not a general existence result of solutions to (\ref{a1}), (\ref{q1}), (\ref{q2}). Notice indeed that the initial data satisfied corresponding to these solutions are all of the form $f(0, m, p)=m^{-3\theta}F_{in}(m)\, \varphi\left(p/m^\theta \right)$. An example of admissible function $a$ is $a(p) := |p|^2$, so that that $\theta = 1/2$. 
\end{rem}

\noindent
\textbf{Proof of Theorem \ref{S4T1}.} 
We have to check that the function $f(t, m, p)$ defined by (\ref{S4E1bis}) solves (\ref{a1}), (\ref{q1}), (\ref{q2}). We start with writing
\bear
\label{S4E3}
&&\frac{\partial f}{\partial t}=m^{-3\theta}\varphi\left(\frac{p}{m^{\theta}} 
\right)\frac{\partial F}{\partial t} \nonumber \\
&&=m^{-3\theta}\varphi\left(\frac{p}{m^{\theta}} 
\right)\left[
\frac{1}{2}\int_0^m\, F(t, m-m') F(t, m')\, a(m-m', m')\,  dm'\right. \nonumber\\
&&\hskip 6cm \left. - \int_0^\infty  F(t, m) F(t, m')\, a(m, m')\,  dm'\right].
\eear
On the one hand, using that 
$$
\int_{\R^3} \varphi (p) \, dp = \FF(\varphi) (0) = e^{-a (0)Ê}Ê= 1,
$$
the last term in (\ref{S4E3}) gives
\bear  \nonumber 
&&m^{-3\theta}\varphi\left(\frac{p}{m^{\theta}}
\right) \int_0^\infty  F(t, m) F(t, m')\, a(m, m')\, dm'=
\\  \nonumber 
&&= m^{-3\theta}\varphi\left(\frac{p}{m^{\theta}} \right)F(t, m) 
\int_0^\infty \, a(m, m')\, F(t, m') 
\int_{\R^3} (m')^{-3\theta} \varphi(\frac{p'}{m'^\theta}) \, dp'
\\  \label{S4E4}
&&=f(t, m, p)\int_0^\infty\int_{\R^3} \, a(m, m')\,  f(t, m', p')\, dp'.
\eear
On the other hand, let us define the function
$$
g(m, p)=m^{-3\, \theta}\varphi(p/m^\theta).
$$
Using the definition of $\varphi$ and the homogeneity of $a$, it satisfies for any $0 < m' < m$
\bean
\hat g (m,\eta) &=& \hat \varphi (m^\theta \, \eta) = \exp (- a(m^\theta \, \eta)) = \exp (- m \, a( \eta)) 
\\
&=& \exp (- m' \, a( \eta)) \, \exp (- (m - m') \, a( \eta)) 
\\
&=& \hat g (m',\eta) \, \hat g (m-m',\eta),
\eean
or coming back to the origin function 
$$
g(m,p) = \int_{\R^3} g(m',p') \, g(m-m',p-p') \, dp'.
$$
Using that identity in the first (gain) term in (\ref{S4E3}), we get 
\bear  \nonumber 
&&m^{-3\theta}\varphi\left(\frac{p}{m^{\theta}}
\right) \int_0^m\, F(t, m-m') F(t, m')\, a(m-m', m')\,  dm'=
\\  \nonumber 
&&=g(m,p) \int_0^m\, F(t, m-m') F(t, m')\, a(m-m', m') \, dm'
\\  \nonumber
&&= \int_{\R^3} \int_0^m F(t, m-m') g(m-m',p-p') \, F(t, m') g(m',p') \, a(m-m', m') \, dm'dp'
\\  \label{S4E5}
&&= \int_{\R^3} \int_0^m f(t, m-m',p-p') \, f(t, m',p') \, a(m-m', m') \, dm'dp'. 
\eear
We conclude that $f$ satisfies \ref{a1}), (\ref{q1}), (\ref{q2}) by gathering (\ref{S4E3}), (\ref{S4E4}) and (\ref{S4E5}). \qed

\medskip

The previous Theorem is useful in order to prove the existence of self similar solutions for some kernels $a(m, m')$ as it is seen in the following corollary.

\begin{cor}
\label{S4T2}
Suppose  that $a$ and $\theta $ are as in Theorem \ref{S4T1}. Assume further that $F$ is a self similar solution of the coagulation equation with coagulation kernel $a(m, m')$. Then the function $f$ defined by  (\ref{S4E1bis}) is a self similar solution of
(\ref{a1}), (\ref{q1}), (\ref{q2}).
\end{cor}

\noindent
\textbf{Proof of Corollary \ref{S4T2}.} The hypothesis on $F$ means that for some functions $\Phi$, $\nu(t)$ and $\mu (t)$ it may be written as:
\bean
F(t, m)=\nu(t)\Phi(\mu(t)\, m).
\eean
Therefore $f$ is a self-similar function since it may be written as 
\bean
f(t, m, p) & = &  m^{-3\theta}\,\nu(t)\Phi(\mu(t)\, m)\, \varphi\left( \frac{p}{m^{\theta}}\right)
\\
& = & \nu(t) \mu(t)^{3\theta}\left(\mu(t)\, m \right)^{-3\theta}\, \Phi(\mu(t)\, m)\varphi\left( \frac{\mu(t)^{\theta}\, p}{(\mu(t)m)^{\theta}}\right)\\
& = & \nu(t) \mu(t)^{3\theta}\, \Psi \left(\mu(t)\, m, \mu(t)^{\theta}\, p \right)
\eean
with $\Psi (M, P) =  M^{-3\theta}\, \Phi(M)\, \varphi\left(P / M^{\theta} \right)$. \qed

\begin{rem}
\label{S4R2}
(i) Self similar solutions for equation (\ref{a1}), (\ref{q1}), (\ref{q2}) had already been obtained in \cite{FM}. They correspond to the case $\theta=1/2$ of the above Corollary.

\smallskip\noindent
(ii) Self similar solutions of the coagulation equation are well known to exist for the cases $a(m, m')=1$, $a(m, m')=m+m'$ and $a(m, m')=m\, m'$. Their existence for several other kernels with homogenetity $\lambda <1$ have been proved in \cite{EMR} and \cite{FL}. In that last case,  these self similar solutions are of the form:
\bear
\label{S4E7}
F(t, m)=t^{-\frac{2}{1-\lambda}}\Phi\left(\frac{m}{t^{\frac{1}{1-\lambda}}} \right).
\eear
We deduce under the assumption of the above Corollary that 
\bear
\label{S4E8}
f(t, m, p)=t^{-\frac{2}{1-\lambda}}m^{-3\theta} \,
\Phi\left(\frac{m}{t^{\frac{1}{1-\lambda}}} \right)\, \varphi\left( \frac{p}{m^{\theta}}\right).
\eear
is a self similar solutions to equation (\ref{a1}), (\ref{q1}), (\ref{q2}) for the same kernel $a(m, m')$. A straightforward calculation yields
\bear
\label{S4E9}
P_k(t)=\int_{\R^d}\int_0^\infty |p|^k\, f(t, m, p)\, dm\, dp=
t^{-\frac{1-k\theta}{1-\lambda}}\int_{\R^d}|P|^k\varphi(P)dP\, \int_0^\infty M^{k\theta}\, \Phi(M)\, dM.
\eear
As a consequence, we have $P_0 \to 0$, $P_1 \to 0$ and more generally $P_k \to 0$ whenever $k < \theta^{-1}$ but  $P_k/P_0 \to \infty$ for any $k > 0$ and $P_k \to \infty$ whenever $k< \theta^{-1}$. The rough physics interpretation is that the total number of particle decreases, the total impulsion of the gas also decreases, but for instance the mean second moment $P_2/P_0$ tends to infinity in the large time asymptotic, which is the opposite behavior with respect to the one discussed in Remark~\ref{S3Rem3}. Here, the behavior is quite similar with the bahavior of the solutions to Smoluchoski equation since the mean impulsion moment $P_k/P_0 \to \infty$ for any $k > 0$. That makes again a difference with the model discussed in Remark~\ref{S3Rem3}. 
\end{rem}


\section{The constant case  $a = 1$}

\setcounter{equation}{0}
\setcounter{theo}{0}

For the sake of simplicity, we restrict our study to the case $d=1$. It is likely that it extends to higher dimension $d \in \NN^*$. 

\begin{theo}
\label{S5T1}
Suppose that the initial data $f_{in}$ is even, regular and good decreasing properties. Then 
(\ref{a1}), (\ref{a2}), (\ref{q1}) (\ref{q2}) has a solution given by:
\bear \label{def:f=FourLapF}
&&f(t, m, p)=\mathcal F^{-1}\left( \mathcal L ^{-1}F\right)(t, m, p)\\ \label{S5E4}
&& F(t, \zeta, \xi)=\frac{H_0^2}{\left(H_0+(t/2) \right)^2\left(\frac{1}{F(0, \zeta, \xi)}-\frac{H_0\, t/2}{H_0+(t/2)} \right)},
\eear
with $H_0 := M_{0,0}(f_{in})^{-1}$ as defined in (\ref{def:Mab}). Furthermore, $f$ satisfies 
\beqn\label{S5:LtimeCvgce1}
t^{5/2}Ê\, f(t, t \, m, \sqrt{t} \, p) \ \rightharpoonup \varphi_\infty(m,p) 
:=  
 \frac{4 \, H^2_0}{\sqrt{2 \, \pi \, \Ac \, \Bc}}
\frac{e^{-\frac{2 \, H^2_0 \, m}{\Ac}} e^{-\frac{\Ac \, p^2}{ 2 \, \Bc \, m}}}{\sqrt{ m}}
\eeqn
in the weak sense of measure $\sigma(M^1(Y), C_c(Y))$ as $ t\to +\infty$, where
\bear
\mathcal A & = & {H^2(0)}\int_0^\infty \int_{\R} m\, f(0, m, p)dp dm,\label{S5E3ter}\\
\mathcal B & = & {H^2(0)\over 2}\int_0^\infty \int_{\R}\, p^2\, f(0, m, p)dp dm. \label{S5E3quatr}
\eear
\end{theo}
\textbf{Proof of Theorem \ref{S5T1}}
We first notice that the equation (\ref{a1}), (\ref{q1}) (\ref{q2}) is now:
\bear
\label{S5E1}
\partial _t f(t, m, p)=\frac{1}{2}\int_{\R^d}\int_0^mf(t, m-m', p-p')\, f(t, m', p')dm'\, dp' \nonumber \\
-f(t, m, p)\, \int_{\R^d}\int_0^\infty f(t, m', p')\, dm'\, dp'.
\eear
This equation may be explicitly solved using Fourier transform with respect to $p\in \R$ and Laplace transform with respect to $m>0$. Of course this needs the transform $F$ of the  function $f$ to be defined. This has then to be checked once the expression of $f$ is obtained. We thus define
\bear
\label{S5E2}
F(t, \zeta, \xi)= \int_0^\infty \int_{\R^d} e^{-m\, \zeta}\, e^{-  i \, p\, \xi}
f(t, m, p)\, dp\, dm.
\eear
We then take formally Fourier and Laplace transforms in (\ref{S5E1}) to obtain the Bernouilli equation:
\bear
\label{S5E3}
&&\partial _t F(t, \zeta, \xi)=\frac{1}{2}F^2(t, \zeta, \xi)-M_0(t)\, F(t, \zeta, \xi)\\
&& M_0(t)=F(t, 0, 0).\label{S5E3bis}
\eear
We first notice, taking $\zeta=\xi=0$ in (\ref{S5E3}), that $M_0(t)$ satisfies $\frac{d}{dt}M_0(t)=-\frac{1}{2}M_0^2(t)$ from where
\bear
\label{S5EMzero}
M_0(t)=\frac{1}{H_0+t/2}.
\eear
Classical ODE integration methods lead that the solution of (\ref{S5E3}) is the function $F(t, \zeta, \xi)$ given by (\ref{S5E4}).
On the one hand, the function $t \mapsto (H_0 \,  t/2)/ (H_0 +  t/2)$ is strictly increasing with limit in infinity equal to $H_0$, so that for any
$\delta \in (0,1)$ there exists $T \in (0,\infty)$
\beqn\label{S5:ineq1}
\forall \, t \in [0,T] \qquad  \frac{H_0 \,  t/2}{H_0 +  t/2} \le H_0 \, (1-\delta),
\eeqn
and on the other hand 
\beqn\label{S5:ineq2}
|F(0,\zeta,\xi)| \le \int_0^\infty \int_\R f(0,m,p) \, dm dp = H_0^{-1}.
\eeqn
Gathering (\ref{S5:ineq1}) and  (\ref{S5:ineq2}) the fraction in the right hand side of (\ref{S5E4}) is well defined for all $t>0$. 
More precisely for any $t \in [0,T]$
\bean
\left|\frac{1}{F(0, \zeta, \xi)}-\frac{H_0\, t/2}{H_0+(t/2)} \right| 
&\ge& \left| \frac{1}{F(0, \zeta, \xi)} \right|  -\frac{H_0\, t/2}{H_0+(t/2)} 
\\
&\ge& \left| {F(0, \zeta, \xi)} \right|^{-1}  - H_0 \, (1-\delta) \ge \delta \,  \left| {F(0, \zeta, \xi)} \right|^{-1}, 
\eean
which implies 
\bear
\label{S5E6}
|F(t, \zeta, \xi)|\le \frac{H_0^2}{\delta \left(H_0+t/2 \right)^2}|F(0, \zeta, \xi)|.
\eear
As a consequence, any ``good'' decay and regularity properties of the initial data $f_{in}$ ensure ``good'' decay and regularity properties of 
$F(0, \zeta, \xi)$. It is then possible to take the inverse Fourier and Laplace transforms of $F(t, \zeta, \xi)$ to define the function $f(t, m, p)$.\\ 
If one is interested in the behaviour of $f(t, m, p)$ as $t\to \infty$ it is a classical argument to consider the rescaled functions $\varphi$ associated to $f$ by the 
relation
\beqn
\varphi(t,M,P) := t^{5/2} \, f(t, t \, M, \sqrt{t} \, P),
\eeqn
so that 
\bear
\label{S5E6}
f(t, m, p)=t^{-5/2}\, \varphi\left(t, \frac{m}{t}, \frac{p}{\sqrt t} \right).
\eear
Taking the Fourier and Laplace transform in both side yields 
\bear
F(t, \zeta, \xi) & = & t^{-1}\Phi (t, t \, \zeta, \sqrt t \, \xi) \label{S5E7}
\eear
with
\bear
\Phi(t, \zeta, \xi) & = & \frac{t H^2_0}
{\left(H_0+(t/2) \right)^2\left(\frac{1}{F(0, \frac{\zeta}{t}, \frac{\xi}{\sqrt t})}-\frac{H_0\, t/2}{H_0+(t/2)} \right)}.\label{S5E8}
\eear
Since we are interested in the long time behaviour of $\Phi(\cdot,\zeta,\xi)$ for all $\zeta$ and $\xi$ fixed , we may write:
\bean
\frac{1}{F(0, \frac{\zeta}{t}, \frac{\xi}{\sqrt t})}-\frac{H_0\, t/2}{H_0+(t/2)} 
& = & 
\frac{1}{F(0, \frac{\zeta}{t}, \frac{\xi}{\sqrt t})}-H_0+\frac{H_0^2}{H_0+(t/2)} 
\eean
and consider the auxiliary function
\bear
\Psi(t, \zeta, \xi) & = & \frac{t H^2_0}
{\left((t/2) \right)^2\left(\frac{1}{F(0, \frac{\zeta}{t}, \frac{\xi}{\sqrt t})}-H_0+\frac{H_0^2}{(t/2)} \right)} \nonumber\\
& = & \frac{4 H^2_0}
{t\left(\frac{1}{F(0, \frac{\zeta}{t}, \frac{\xi}{\sqrt t})}-H_0+\frac{2\,H_0^2}{t} \right)}.\label{S5E9}
\eear
We perform the following expansion up to the order $o(1/t)$: 
\bear
\label{S5E10}
\frac{1}{F(0, \frac{\zeta}{t}, \frac{\xi}{\sqrt t})}-H_0=\frac{\zeta}{t} \frac{\partial F^{-1}}{\partial \zeta}(0, 0, 0)+\frac{\xi}{\sqrt  t} \frac{\partial F^{-1}}{\partial \xi}(0, 0, 0)+\frac{1}{2}\frac{\xi^2}{  t} \frac{\partial ^2 F^{-1}}{\partial \xi ^2}(0, 0, 0)+o\left(\frac{1}{t} \right)\!. \quad
\eear
Since by hypothesis $f$ is even with respect to $p$,  we have
\bean
\frac{\partial F}{\partial \xi}(0, 0, 0)=-i\int_0^\infty\int_{\R} f_{in}(m, p) \, p\, dp\, dm=0
\eean
and  then:
\beqn\label{S5:eqFxi}
\frac{\partial F^{-1}}{\partial \xi}(0, 0, 0)=-\frac{1}{F(0,0,0)^2}\frac{\partial F}{\partial \xi}(0, 0, 0) =0.
\eeqn
We also have 
\bean
\frac{\partial^2 F}{\partial \xi^2}(0, 0, 0) = -\int_0^\infty \int_{\R}\, p^2\, f(0, m, p)dp dm,
\eean
which with the help of (\ref{S5:eqFxi}) implies 
\bear\nonumber 
\frac{\partial^2 F^{-1}}{\partial \xi^2} (0,0,0) & = & -F^{-2}(0,0,0)\frac{\partial^2 F}{\partial \xi^2} (0,0,0)+F^{-3}(0,0,0)\left( \frac{\partial F}{\partial \xi} (0,0,0) \right)^2
\\ \label{S5:ident2}
& = & {H^2(0)} {\int_0^\infty \int_{\R}\, p^2\, f(0, m, p)dp dm} \,\,\, = 2\, \mathcal B.
\eear
Similarly, we compute 
\bean
\frac{\partial F}{\partial \zeta}(0, 0, 0) & = & -\int_0^\infty \int_{\R}m\, f(0, m, p)dp dm,
\eean
which implies
\beqn\label{S5:ident3}
\frac{\partial F^{-1}}{\partial \zeta}(0, 0, 0)=  -\frac{1}{F^2(0, 0, 0)}\frac{\partial F}{\partial \zeta}(0, 0, 0)= \mathcal A.
\eeqn
Thanks to (\ref{S5E10}), (\ref{S5:eqFxi}), (\ref{S5:ident2}) and (\ref{S5:ident3}), we deduce that (\ref{S5E9}) reads now:
\bean
\Psi(t, \zeta, \xi)=\frac{4 H^2(0)}
{\left( \zeta \, \Ac + \xi^2 \, \Bc +2 H^2(0) +o(1)\right) }
\eean
from where
\bear
\label{S5E12}
\lim_{t\to +\infty} \Phi(t, \zeta, \xi)=\lim_{t\to +\infty} \Psi(t, \zeta, \xi)=
\frac{4 H^2(0)}
{ \mathcal A\, \zeta +\mathcal B\, \xi^2+2 H^2(0) } =: \Psi_\infty(\zeta,\xi).
\eear
In order to come back to the original variables, we recall that from standard integral calculus for any $\Cc, \Dc > 0$
\bean
\frac{1}{(2\pi)^{1/2}} 
\int_0^\infty \! \int_{\R} e^{-m\, \zeta} \, e^{-i\, p\, \xi} \,\frac{e^{- \Cc \, m } e^{-\frac{|p|^2}{2 \, \Dc \, m}}} {\sqrt {\Dc \, m}}\, dp\, dm
=\frac{1}{ \zeta+\Dc \, \xi^2+\Cc}, 
\eean
from where choosing $\Cc:= 2\, H^2_0/\Ac$ and $\Dc := \Bc/\Ac$,  we obtain
$$
({\cal F}^{-1} {\cal L}^{-1}) \left(  \Psi_\infty \right)= {4 \, H^2_0 \over (2\pi)^{1/2} \, \Ac}  \, \frac{e^{- \Cc \, m } e^{-\frac{|p|^2}{2 \, \Dc \, m}}} {\sqrt {\Dc \, m}}\
= \varphi_\infty (m,p)
$$
as defined in (\ref{S5:LtimeCvgce1}). Finally, (\ref{S5E12}) implies that $ \varphi(t,.) \ \wto \ \varphi_\infty$ in the weak sens of measure, which is nothing but (\ref{S5:LtimeCvgce1}). 
\qed

\medskip
The previous Theorem shows the convergence of some of the solutions of (\ref{a1})-(\ref{q2}) to a function which is a self similar solution of the equation (\ref{a1}), (\ref{q1}), (\ref{q2}), i.e. a solution of the form
\bear
\label{S5E12}
f(t, m, p)=t^{-\alpha}\varphi(t^{-1}\, m, t^{-\beta}\, p)
\eear
for some function $\varphi$. The numbers $\alpha$ and $\beta$ define the scaling of the self similar solutions. In the Theorem \ref{S5T1} we have $\alpha=5/2$ and $\beta =1/2$.
It turns out that equation (\ref{a1}), (\ref{q1}), (\ref{q2}) has more than one self similar solution with the same scaling as it is shown in the next Theorem.

\begin{theo}
\label{S5T2}
Let $\Phi\in C^1(\R^d)$ is such that 
$$
g(y, x)=\mathcal F _{\xi}^{-1}\mathcal L _{\zeta}^{-1}\left( \frac{2}{2\, \zeta \Phi \left(\frac{\xi^2}{\zeta}\right)+1}\right)
$$
is well defined for $x \in \R^d$ and $y>0$. Then
\beqn\label{S5:SSg}
t^{-\frac{5}{2}}g\left(\frac{m}{t},\, \frac{p}{\sqrt t} \right).
\eeqn
is a self similar solution to (\ref{a1}), (\ref{q1}), (\ref{q2}). 
\end{theo}

\noindent
\textbf{Proof of Theorem \ref{S5T2}.} We look after self similar solutions of the form (\ref{S5:SSg}).
The function $g$ must then solve:
\bear
\label{S5E13}
-\frac{5}{2}\, g-y\partial _y g-\frac{1}{2}x\partial_x g=
\frac{1}{2}\int_{\R}\int_0^y g(y-y', x-x')\, g(y', y')dy'\, dx'-\nonumber \\ 
-g \int_{\R}\int_0^\infty g(y',x')dy'\, dx'. 
\eear
We integrate this equation with respect to $x$ and $y$ and obtain
\bear
\label{S5E12bis}
\int_{\R}\int_0^\infty g(y',x')dy'\, dx'=2.
\eear
We now Fourier transform with respect to $x$ and Laplace transform with respect to $y$:
\bean
\zeta \partial_\zeta \widehat g+\frac{1}{2}\xi \partial_\xi \widehat g=\frac{1}{2}\widehat g^2-\widehat g.
\eean
We divide by $\widehat g^2$ and define $G=1/\widehat g$:
\bean
\zeta\partial_\zeta G+\xi \partial_\xi G=G-\frac{1}{2}.
\eean
The function $G$ may then be any function of the form:
\bean
G(\zeta, \xi)=\zeta \Phi \left(\frac{\xi^2}{\zeta}\right)+\frac{1}{2}
\eean
for any arbitrary derivable function $\Phi$. Therefore
\bear
\label{S5E14}
\widehat g(\zeta, \xi)=\frac{2}{2\, \zeta \Phi \left(\frac{\xi^2}{\zeta}\right)+1},
\eear
with, due to (\ref{S5E12bis}):
\bean
\lim_{\zeta\to 0,\, \xi \to 0}\frac{2}{2\, \zeta \Phi \left(\frac{\xi^2}{\zeta}\right)+1}=2\Longleftrightarrow
\lim_{\zeta\to 0,\, \xi \to 0}\zeta \Phi \left(\frac{\xi^2}{\zeta}\right)=0.
\eean
If we want to define the function $g$ from (\ref{S5E14}) the function $\Phi$ must be such that $\widehat g$ has an inverse Fourier and Laplace transform. \qed

\begin{rem}
If $\Phi(z)=z+1$, 
\bean
\hat g(\zeta, \xi)=\frac{2}{2\, \zeta \left(\frac{\xi^2}{\zeta}+1\right)+1}
=\frac{2}{2\,  (\xi^2+\zeta)+1}\\
={\cal L} \left({\cal F}\left(\frac{e^{-\frac{y}{2}}\, e^{-\frac{x^2}{4 y}}}
{\sqrt 2 \sqrt y} \right)\right).
\eean
This is the profile of the self similar solution which appears in Theorem \ref{S5T1}.
It is easy to obtain particular solutions $g$, some of them are explicit others are not.
If, for example, $\Phi\equiv 1$ then $g(y, x)= e^{-y^2}\delta_{x=0}$. Another explicit example is for  $\Phi(z)=z$ which gives
$g(x, y)=\sqrt \pi \delta_{y=0}e^{-\frac{|x|}{\sqrt 2}}$. On the other hand, if we take $\Phi(z)=\sqrt z$, the inverse Laplace transform, let us call it $h(y, \xi)$, is still explicit:
\bear
\label{S5E15}
h(y, \xi)=\mathcal L^{-1}_{\zeta}\left(\frac{2}{2\sqrt{\zeta\, \xi^2}}+1 \right)
=\frac{\frac{\sqrt{\xi^2}}{\sqrt \pi\, \sqrt y}-e^{\frac{y}{\xi^2}}Erfc\left( \sqrt{\frac{y}{\xi^2}}\right)}{4\xi^2}.
\eear
\end{rem}
It remains to check that $h(y, \cdot)$
has an inverse Fourier transform with respect to the variable $\xi$. It is easily checked that, for all $y>0$ fixed:
\bean
&&h(y, \xi)=\mathcal O \left(\frac{\xi}{y^{3/2}} \right), \,\,\,\hbox{as}\,\,\xi\to 0\\
&&h(y, \xi)=
\frac{\frac{1}{\sqrt \pi}\sqrt{\frac{\xi^2}{y}}-1}{4\xi^2}+\mathcal O\left(\frac{y}{\xi^2} \right)
\,\,\,\hbox{as}\,\,|\xi|\to +\infty.
\eean
This function is then in $L^2(\R)$ with respect to the $\xi$ variable and has then an inverse Fourier transform with respect to $\xi$ which is $g(y, x)$:
\bean
g(y, x)=\mathcal F^{-1}_{\xi}(h(y, \cdot))(x).
\eean
 Moreover, for all $y>0$, $g(y,\cdot)\in L^2(\R)$ and the convolution of $g(y,\cdot)$ with itself is well defined
$$
\mathcal F \left( g(y-y',\cdot)*g(y',\cdot)\right)(\xi) = h(y-y', \xi)\, h(y', \xi)
$$
and
$$
\int_0^y \mathcal F \left( g(y-y',\cdot)*g(y',\cdot)\right)(\xi)dy = 
\int_0^y h(y-y', \xi)\, h(y', \xi)\, dy.
$$
Therefore,
\bean
\int_{\R}\left|\mathcal F \left( g(y-y',\cdot)*g(y',\cdot)\right)(\xi)\right|
d\xi & \le & 
\int_0^y \int_{\R}\left|h(y-y', \xi)\, h(y', \xi)\right| d\xi\, dy = \sum_{k=1}^6 I_k,
\eean
with
\bean
&&I_1:=\int_0^{y/2} \int_{|\xi|\le {y'}^{1/2} \le {(y-y')}^{1/2}}\left|h(y-y', \xi)\, h(y', \xi)\right| d\xi\, dy',\\
&&I_2:=\int_0^{y/2} \int_{{y'}^{1/2}\le|\xi|\le {(y-y')}^{1/2}}\left|h(y-y', \xi)\, h(y', \xi)\right| d\xi\, dy',\\
&&I_3:=\int_0^{y/2} \int_{{y'}^{1/2}\le {(y-y')}^{1/2}\le |\xi|}\left|h(y-y', \xi)\, h(y', \xi)\right| d\xi\, dy',\\
&&I_4:=\int_{y/2}^y \int_{|\xi|\le {(y-y')}^{1/2}\le {y'}^{1/2} }\left|h(y-y', \xi)\, h(y', \xi)\right| d\xi\, dy',\\
&&I_5:=\int_{y/2}^y \int_{{(y-y')}^{1/2}\le|\xi|\le {y'}^{1/2}}\left|h(y-y', \xi)\, h(y', \xi)\right| d\xi\, dy',\\
&&I_6:=\int_{y/2}^y \int_{{(y-y')}^{1/2}\le {y'}^{1/2}\le |\xi|}\left|h(y-y', \xi)\, h(y', \xi)\right| d\xi\, dy'.
\eean
We must verify that each term is finite. Indeed, we have 
\bean
I_1 & \le & C\, \int_0^{y/2} \int_{|\xi|\le {y'}^{1/2}\le (y-y')^{1/2}}\frac{\xi^2}{{y'}^{3/2}{(y-y')}^{3/2}} d\xi\, dy'\\
& \le & C\, \int_0^{y/2} \frac{min\{{y'}^{3/2},\,(y-y')^{3/2}\}}{{y'}^{3/2}{(y-y')}^{3/2}}dy'<\infty ; \\
I_2 & \le & C \int_0^{y/2} \int_{{y'}^{1/2}\le|\xi|\le {(y-y')}^{1/2}}\frac{|\xi|}{(y-y')^{3/2}}\,\left( \frac{1}{\sqrt {y'}\, |\xi|}+\frac{1}{\xi^2}+\mathcal O\left(\frac{y'}{\xi^2} \right)\right) d\xi\, dy'\\
& \le & \frac{C}{y^{3/2}}\int_0^{y/2} \int_{{y'}^{1/2}\le|\xi|\le {(y-y')}^{1/2}}
\left(\frac{1}{\sqrt {y'}}+\frac{1}{|\xi|} +1\right)d\xi dy\\
& \le & \frac{C}{y^{3/2}}\int_0^{y/2} \int_{{y'}^{1/2}\le|\xi|\le {(y-y')}^{1/2}}
\left(\frac{2}{\sqrt {y'}}+1\right)d\xi dy<\infty ; \\
I_3 & \le & C \int_0^{y/2} \int_{{y'}^{1/2}\le {(y-y')}^{1/2}\le |\xi|} 
\left( \frac{1}{\sqrt {y'}\, |\xi|}+\frac{1}{\xi^2}+\mathcal O\left(\frac{y'}{\xi^2} \right)\right)\times \\
&&\times \left( \frac{1}{\sqrt {y-y'}\, |\xi|}+\frac{1}{\xi^2}+\mathcal O\left(\frac{y-y'}{\xi^2} \right)\right)d\xi\, dy'\\
& \le & C\int_0^{y/2} \int_{{y'}^{1/2}\le {(y-y')}^{1/2}\le |\xi|}
\left(\frac{1}{\sqrt {y'}\sqrt {y-y'}\, |\xi|^2} 
+\frac{1}{|\xi|^3}\left(\frac{1}{\sqrt {y-y'}}+\frac{1}{\sqrt {y'}} \right)+\right. \\
&&\left.+\frac{1}{\xi^4}+\mathcal O \left(\frac{y'}{\sqrt{y-y'}|\xi|^3} \right)+\mathcal O \left(\frac{y-y'}{\sqrt{y'}|\xi|^3} \right)
+\mathcal O \left(\frac{y+y^2}{|\xi|^4} \right)\right)d\xi dy\\
& \le & \frac{C}{\sqrt y}\int_0^{y/2} \frac{1}{\sqrt {y'}}\int_{{y'}^{1/2}\le {(y-y')}^{1/2}\le |\xi|}\frac{d\xi}{|\xi|^2}\, dy+\\
&&+C\int_0^{y/2}\left(\frac{1}{\sqrt {y}}+\frac{1}{\sqrt {y'}} \right) \int_{{y'}^{1/2}\le {(y-y')}^{1/2}\le |\xi|}\frac{d\xi}{|\xi|^3}\, dy+\\
&&+C\int_0^{y/2} \int_{{y'}^{1/2}\le {(y-y')}^{1/2}\le |\xi|} \frac{d\xi}{|\xi|^4}\, dy\\
& \le & \frac{C}{\sqrt y}\int_0^{y/2} \frac{dy}{\sqrt {y'}\, \sqrt {(y-y')}}+
C\int_0^{y/2}\left(\frac{1}{\sqrt {y}}+\frac{1}{\sqrt {y'}} \right)
\frac{dy}{y-y'}+C\int_0^{y/2}\frac{dy}{(y-y')^2}.
\eean
Similar estimates show that the integrals $I_4, I_5$ and $I_6$ converge. The function $\int_0^y h(y-y', \xi)\, h(y', \xi)\, dy$ is then in $L^1(\R)$ and has then an inverse Fourier transform which is 
\bean
\int_0^y (g(y-y', \cdot)*g(y', \cdot))(\xi)\, dy.
\eean

\end{document}